\documentclass[11pt]{article}
\usepackage{amsfonts,amsmath,amssymb,accents}
\usepackage{hyperref}
\usepackage{graphicx}

\newcommand{\figfilebw}[1]{bitfigures/#1}

\newcommand{\figfile}[1]{figures/#1}

\newcommand{\detail}[1]{\par\noi{\bf [Proof detail\ }{#1}
\hfill{\bf ]}\par\noi\hspace{-4pt}}
\renewcommand{\detail}[1]{}

\newcommand{\dis}{\displaystyle}

\newcommand{\noi}{\noindent}

\newtheorem{theorem}{Theorem}
\newtheorem{proposition}[theorem]{Proposition}
\newtheorem{corollary}[theorem]{Corollary}
\newtheorem{conjecture}[theorem]{Conjecture}
\newtheorem{lemma}[theorem]{Lemma}

\newtheorem{remark}[theorem]{Remark}

\newcommand{\bt}{\begin{theorem}}
\newcommand{\et}{\end{theorem}}
\newcommand{\bl}{\begin{lemma}}
\newcommand{\el}{\end{lemma}}
\newcommand{\bp}{\begin{proposition}}
\newcommand{\ep}{\end{proposition}}
\newcommand{\bcor}{\begin{corollary}}
\newcommand{\ecor}{\end{corollary}}
\newcommand{\br}{\begin{remark}\rm}
\newcommand{\er}{\end{remark}}
\newcommand{\bcon}{\begin{conjecture}}
\newcommand{\econ}{\end{conjecture}}

\newcommand{\be}{\begin{equation}}
\newcommand{\ee}{\end{equation}}
\newcommand{\ba}{\begin{array}}
\newcommand{\ea}{\end{array}}
\newcommand{\bc}{\be\begin{array}{r@{\,}c@{\,}l}}
\newcommand{\ec}{\end{array}\ee}

\newcommand{\al}{\alpha}
\newcommand{\bet}{\beta}

\newcommand{\de}{\delta}
\newcommand{\De}{\Delta}

\newcommand{\tet}{\theta}

\newcommand{\R}{{\mathbb R}}

\newcommand{\Z}{{\mathbb Z}}

\renewcommand{\P}{{\mathbb P}}
\newcommand{\E}{{\mathbb E}}

\newcommand{\down}{\downarrow}
\newcommand{\sub}{\subset}

\newcommand{\un}{\underline}

\newcommand{\rvec}[1]{\accentset{\rightarrow}{#1}}
\newcommand{\lvec}[1]{\accentset{\leftarrow}{#1}}

\newcommand{\ffrac}[2]{{\textstyle\frac{{#1}}{{#2}}}}
\newcommand{\dif}[1]{\ffrac{\partial}{\partial{#1}}}

\newcommand{\di}{\mathrm{d}}

\begin{document}

\makeatletter\@addtoreset{equation}{section}
\makeatother\def\theequation{\thesection.\arabic{equation}} 

\renewcommand{\labelenumi}{{(\roman{enumi})}}

\title{\vspace{-3cm}Numerical analysis of the rebellious voter model}
\author{Jan M.~Swart${}^\ast$ \and Karel Vrbensk\'y
\footnote{swart@utia.cas.cz and vrbensky@utia.cas.cz,
Institute of Information Theory and Automation of the ASCR (\' UTIA),
Pod vod\'arenskou v\v e\v z\' i 4, 18208 Praha 8,
Czech Republic. URL http://staff.utia.cas.cz/swart/}}

\maketitle\vspace{-30pt}

\begin{abstract}\noi
The rebellious voter model, introduced by Sturm and Swart (2008), is a
variation of the standard, one-dimensional voter model, in which types that
are locally in the minority have an advantage. It is related, both through
duality and through the evolution of its interfaces, to a system of branching
annihilating random walks that is believed to belong to the
`parity-conservation' universality class. This paper presents numerical data
for the rebellious voter model and for a closely related one-sided version of
the model. Both models appear to exhibit a phase transition between
noncoexistence and coexistence as the advantage for minority types is
increased. For the one-sided model (but not for the original, two-sided
rebellious voter model), it appears that the critical point is exactly a half
and two important functions of the process are given by simple, explicit
formulas, a fact for which we have no explanation.
\end{abstract}

\noi
{\it Acknowledgements.} Work sponsored by GA\v CR grant 201/06/1323 and
M\v{S}MT \v{C}R grant DAR 1M0572. We thank Ivan Dornic, as well as two
referees, for many useful comments, especially on the physical side of the
problem.

{\small\setlength{\parskip}{-2pt}\tableofcontents}
\newpage

\section{Introduction}

The study of parity preserving branching-annihilating particle systems has a
long history in the physics literature. Their main interest lies in the fact
that they provide a rare example of systems that exhibit a phase transition of
survival/extinction type (the physics literature usually speaks of a
transition between an `active' and an `absorbing' phase) that does not belong
to the directed percolation (DP) universality class.

Indeed, the phase transition of these models appears to belong to a different
`parity conservation' (PC) universality class, characterized by its own
critical exponents. The first example of a phase transition belonging to this
universality class was described in \cite{GKT84}; this concerned a
one-dimensional model for interfaces, motivated by solid state
physics. Numerical simulations for parity preserving branching-annihilating
particle systems on $\Z^d$ in dimensions $d=1,2,3$ and also on Sierpinski
gaskets can be found in \cite{TT92}. Further simulations for a one-dimensional
model were made in \cite{Jen94}. This paper also contains the conjecture that
the order parameter critical exponent $\bet$ equals $13/14$; in another paper
the author has conjectured that $\bet=1$. The latter value was also found in
\cite{IT98}, while more recent literature seems to agree on a value of
$\bet\approx 0.92$ or slightly higher, see e.g.\ \cite{Hin00,OS05}. Our
present work is relevant to the question of determining the value of $\bet$
since our proposed explicit formula for the survival probability $\rho$ of the
one-sided rebellious voter model, if it is correct, implies that $\bet=1$.

Field-theoretic (renormalization group) work on the PC universality class was
carried out in \cite{CT98}, who found two critical dimensions. The (usual)
upper critical dimension is $d_{\rm c}=2$, above which one should find
mean-field exponents. There is, however, a second critical dimension $d'_{\rm
  c}\approx 4/3$ such that only in dimensions below $d'_{\rm c}$ there is a
nontrivial absorbing phase. From the field-theoretical perspective, the
existence of a second critical dimension leads to considerable difficulties,
to which the authors of \cite{CCDDM05} reacted by considering non-perturbative
techniques.

As usual, mathematical theory lags far behind theoretical physics and occupies
itself with more basic and partially different questions. In particular,
renormalization group techniques are generally nonrigorous and the existence
of critical exponents or their universality below the upper critical dimension
are from the mathematical perspective unproved conjectures. Nevertheless, in
recent years, some rigorous mathematical work has been carried out on
models (supposedly) belonging to the PC universality class. Sudbury
\cite{Sud90} studied a `double branching annihilation random walk' (DBARW),
which is a one-dimensional model for parity preserving branching-annihilating
particles. The DBARW that can be exactly solved, in a certain sense, but
unfortunately the model is in the absorbing phase for all nontrivial values of
its parameter, hence there is no `true' phase transition. (Nevertheless, its
solution has been used in papers such as \cite{CT98} to draw conclusions that
presumably hold for the whole PC universality class.) No solvable models seem
to be known so far that exhibit a true phase transition. (Below, we will
present a possible candidate for such a model, the one-sided rebellious voter
model.)

Strong impetus to the mathematical study of models from the PC universality
class came from the paper of Neuhauser and Pacala \cite{NP99}, which proposes
a variation on the $d$-dimensional range-$N$ voter model motivated by
applications in population biology. Their model is governed by two parameters
$\al_0,\al_1$. Here $\al_i\leq 1$ models the {\em interspecific competition
  rate}, i.e., the death rate of organisms of type $i$ due to competition with
organisms of the other type, while the {\em intraspecific competition rate} is
supposed to be one. In the symmetric case $\al_0=\al_1=\al$, their model is
dual to a system of parity preserving branching-annihilating random walks,
with branching rate $1-\al$. For $d=1,N=1$, up to a rescaling of time, this
dual model is the DBARW; results for the latter imply that the model of
Neuhauser and Pacala with $d=1,N=1$ exhibits noncoexistence for all
$\al>0$. For $N\geq 2$ and $d=1$, their model is supposed to exhibit a
nontrivial phase transition between coexistence and noncoexistence. For $d=2$,
their model is supposed to exhibit noncoexistence only for $\al=1$ while in
dimensions $d\geq 3$ noncoexistence holds for all $\al$, in line with the
predictions made in \cite{CT98}.

Much of this lacks rigorous proof. First of all, it is not rigorously known
that coexistence is a decreasing property in $\al$, or, in terms of the dual
model, that increasing the branching rate leads to more particles in the
system. Obvious as this may seem, finding a rigorous proof seems
hard. Second, there exists no rigorous proof of noncoexistence (or
equivalently, the existence of an absorbing phase for the dual model) except
in the trivial cases $\al=1$ (pure voter) or $d=1,N=1$ (related to DBARW). On
the other hand, there exist a number of rigorous results proving coexistence
(or equivalently, the existence of an active phase for the dual
model). Neuhauser and Pacala proved that their model exhibits coexistence for
$\al$ sufficiently close to zero for all values of $d$ and $N$ except
$d=1,N=1$. In \cite{CP06}, it has been proved that the model of Neuhauser and
Pacala exhibits coexistence in dimensions $d\geq 3$ for $\al$ sufficiently
close to {\em one} (the case of intermediate $\al$ is still open due to the
lack of a proof of monotonicity in $\al$). The analogue result in $d=2$ is
proved in \cite{CMP09}.

Some similar models were introduced in \cite{BEM06}, who studied parity
preserving branching-annihilating random walks which allow for more particles
at one site, and \cite{SS08hv}, who introduced the rebellious voter model that
will be studied in the present paper. A simple numerical simulation of that
model, run for illustrative purposes, seemed to suggest that the critical
point of this model is exactly a half, which motivated the present work. In
fact, our results show that the critical value for this particular model is
not $\frac{1}{2}$, but a closely related, one-sided version of the model seems
to have its critical point exactly at $\al=\frac{1}{2}$ and may in some sense
be exactly solvable.

\section{Set-up and basic properties}

\subsection{Definition of the models}

Let
\be
\{0,1\}^\Z:=\big\{x=(x(i))_{i\in\Z}:x(i)\in\{0,1\}\ \forall i\in\Z\big\}
\ee
be the space whose elements $x$ are infinite arrays of zeroes and ones,
situated on the integer lattice $\Z$. We will be interested in Markov
processes taking values in $\{0,1\}^\Z$. We will mainly focus on the {\em
  rebellious voter model} introduced in \cite{SS08hv}, and a variation on this
model, which we call the {\em one-sided rebellious voter model}. For
expository reasons, we will also shortly discuss two other models, the {\em
  disagreement} and {\em swapping voter models}.

In the one-sided rebellious voter model, initially, each site $i\in\Z$ has a
type $x(i)\in\{0,1\}$ assigned to it, and these types are updated according to
the following stochastic dynamics. Independently for each site $i\in\Z$, one
constructs times $0<\tau_1(i)<\tau_2(i)<\cdots$ according to a Poisson
process with rate $1$, i.e.,
$\tau_1(i)-0,\tau_2(i)-\tau_1(i),\tau_3(i)-\tau_2(i),\ldots$ are
i.i.d.\ exponentially distributed with mean $1$. At each time $t=\tau_k(i)$,
with $k=1,2,\ldots$, the type of site $i$ is updated according to the
following rules. With probability $\al$, the site $i$ looks at the site $i-1$
on its left, and if the type of site $i-1$ is different from the type of $i$,
then site $i$ changes its type. With the remaining probability $1-\al$, the
site $i$ looks at the two sites $i-2$ and $i-1$ on its left, and if the type
of site $i-2$ is different from the type of $i-1$, then site $i$ changes its
type. If we let $X_t(i)$ denote the type of site $i$ at time $t$, then
$(X_t)_{t\geq 0}$ is a continuous-time Markov process taking values in
$\{0,1\}^\Z$, which we call the {\em one-sided rebellious voter
  model}. Another way of describing its dynamics is to say that in this model,
for any $i\in\Z$, the process $X$ makes the transition
\be\ba{l}\label{osRV}
x\mapsto x^{\{i\}}\quad\mbox{with rate}\\[5pt]
\dis\quad\al 1_{\{x(i-1)\neq x(i)\}}+(1-\al)1_{\{x(i-2)\neq x(i-1)\}},
\ec
where $1_A$ denotes the indicator function of an event $A$ and for any
$x\in\{0,1\}^\Z$ and $\De\sub\Z$, we let
\be\label{xDe}
x^\De(j):=\left\{\ba{ll}
1-x(j)\quad&\mbox{if }j\in\De,\\
x(j)\quad&\mbox{if }j\not\in\De
\ea\right.\ee
denote the configuration obtained from $x$ by changing the types of all sites
in $\De$.

The original {\em rebellious voter model} as introduced in \cite{SS08hv} is
similar to the process described above, except that sites do not only look to
the left, but also to the right when deciding whether to update their
type. More precisely, in this case, at each time $\tau_k(i)$ as defined above,
the site $i$ decides, with probability $1/2$ each, to look either to the left or
to the right. If the site looks to the left, then its type is updated as
before. If the site $i$ looks to the right, then with probability $\al$, it
looks at the site $i+1$ on its right, and if the types of $i$ and $i+1$
are different, then $i$ changes its type. With the remaining probability
$1-\al$, the site $i$ looks at the two sites $i+1$ and $i+2$ on its right, and
if the types of $i+1$ and $i+2$ different from each other, then $i$ changes
its type. Another way of describing these dynamics is to say that in this
model, for any $i\in\Z$, the process makes the transition
\be\ba{l}\label{RV}
x\mapsto x^{\{i\}}\quad\mbox{with rate}\\[5pt]
\dis\quad\ffrac{1}{2}\al\big(1_{\{x(i-1)\neq x(i)\}}
+1_{\{x(i)\neq x(i+1)\}}\big)\\[5pt]
\dis\quad+\ffrac{1}{2}(1-\al)\big(1_{\{x(i-2)\neq x(i-1)\}}
+1_{\{x(i+1)\neq x(i+2)\}}\big).
\ec
We note that our description of the rebellious voter model differs a factor
$1/2$ in the speed of time in comparison to the original definition in
\cite{SS08hv}.

A third and fourth model, that we also shortly discuss in this introduction, are
the {\em disagreement} and {\em swapping voter models}. The disagreement voter
model evolves as
\be\ba{l}\label{dis}
x\mapsto x^{\{i\}}\quad\mbox{with rate}\\[5pt]
\dis\quad\al\big(1_{\{x(i-1)\neq x(i)\}}+1_{\{x(i)\neq x(i+1)\}}\big)
+(1-\al)1_{\{x(i-1)\neq x(i+1)\}},
\ec
while the swapping voter model evolves as
\be\ba{lll}\label{swap}
x\mapsto x^{\{i\}}\quad&\mbox{with rate}\quad
&\dis\al\big(1_{\{x(i-1)\neq x(i)\}}+1_{\{x(i)\neq x(i+1)\}}\big),\\[5pt]
x\mapsto x^{\{i,i+1\}}\quad&\mbox{with rate}\quad
&\dis(1-\al)1_{\{x(i)\neq x(i+1)\}}.
\ec

In case of the rebellious voter model and its one-sided counterpart, the
parameter $0\leq\al\leq 1$ can be interpreted as the (interspecific) {\em
  competition parameter}. For $\al=1$, the model is a standard,
one-dimensional voter model, as first introduced in \cite{CS73,HL75}. For
$\al<1$, one can check that the dynamics give an advantage to types that are
locally in the minority. As a consequence, at least on a heuristic level,
decreasing $\al$ should make it harder for any type to die out. The rebellious
voter model has been introduced in \cite{SS08hv} in an attempt to model the
distribution of two closely related species where competition between
organisms belonging to different species is less strong than competetion
between organisms of the same species.

\subsection{Interface and dual models}

The models introduced so far are in two ways related to parity preserving
branching annihilating particle systems. First, if $X$ is any of the models
defined in (\ref{osRV}) and (\ref{RV})--(\ref{swap}), then setting
\be\label{inter}
Y_t(i):=1_{\{X_t(i)\neq X_t(i+1)\}}\qquad(t\geq 0,\ i\in\Z)
\ee
defines a Markov process taking values in $\{0,1\}^\Z$ that we call the {\em
  interface model} associated with $X$. (The sites $i$ such that $X_t(i)\neq
X_t(i+1)$ are called interfaces, or also {\em kinks} or {\em domain walls}.)
Indeed, in case of the one-sided rebellious voter model, it is not hard to see
that $Y$ jumps as
\be\label{onesidY}
y\mapsto y^{\{i,i+1\}}\quad\mbox{with rate}
\quad\al1_{\{y(i)=1\}}+(1-\al)1_{\{y(i-1)=1\}},
\ee
where we use notation defined in (\ref{xDe}). For the two-sided process, the
dynamics of $Y$ are given by
\be\ba{l}\label{twosidY}
y\mapsto y^{\{i,i+1\}}\quad\mbox{with rate}\\[5pt]
\dis\quad\ffrac{1}{2}\al\big(1_{\{y(i)=1\}}+1_{\{y(i+1)=1\}}\big)
+\ffrac{1}{2}(1-\al)\big(1_{\{y(i-1)=1\}}+1_{\{y(i+2)=1\}}\big)
\ec
(see \cite[Section~1.2]{SS08hv}), while for the disagreement voter model we
get the `swapping annihilating random walk' (SARW)
\be\label{SARW}
y\mapsto y^{\{i,i+1\}}\quad\mbox{with rate}
\quad\al\big(1_{\{y(i)=1\}}+1_{\{y(i+1)=1\}}\big)+(1-\al)1_{\{y(i)\neq y(i+1)\}},
\ee
and for swapping voter model we get
\be\ba{lll}\label{DBARW}
\dis y\mapsto y^{\{i,i+1\}}\quad&\mbox{with rate}
\quad&\dis\al\big(1_{\{y(i)=1\}}+1_{\{y(i+1)=1\}}\big)\\[5pt]
\dis y\mapsto y^{\{i-1,i+1\}}\quad&\mbox{with rate}\quad&\dis(1-\al)1_{\{y(i)=1\}}.
\ec
Following \cite{Sud90}, we call this latter model the double branching
annihilating random walk (DBARW); very similar models (usually in discrete
time) have been called BARW2 or BAW in the physics literature (see, e.g.,
\cite{TT92}).

In each of these interface models, since $Y$ always flips the types of two
sites at once, it is easy to see that $Y$ {\em preserves parity}, i.e., if $Y$
is started in an initial state $Y_0$ which contains a finite, even
(resp.\ odd) number of ones, then $Y_t$ contains an even (resp.\ odd) number
of ones for each $t\geq 0$. In particular, if $Y$ is started with an odd
number of ones, then the ones can never completely die out.

There is a second relation between our models and parity preserving
branching-annihilating particle systems, namely, through duality.
Recall that two Markov processes $X$ and $Y$ with state spaces $S_X$ and $S_Y$
are {\em dual} to each other if there exists a function $\psi$ defined on
$S_X\times S_Y$, called {\em duality function}, such that
\be\label{Markdual}
\E^x[\psi(X_t,y)]=\E^y[\psi(x,Y_t)]\qquad(t\geq 0),
\ee
where $\E^x$ (resp.\ $\E^y$) denotes expectation with respect to the law of the
process $X$ (resp.\ $Y$) started in $X_0=x$ (resp.\ $Y_0=y$). A necessary and
sufficient condition for (\ref{Markdual}) is that
\be\label{Gdual}
G_X\psi(\,\cdot\,,y)(x)=G_Y\psi(x,\,\cdot\,)(y)\qquad(x\in S_X,\ y\in S_Y),
\ee
where $G_X$ (resp.\ $G_Y$) denotes the generator of $X$ (resp.\ $Y$). In
typical applications, one considers $\psi$ such that the functions
$\big(\psi(\,\cdot\,,y)\big)_{y\in S_Y}$ and $\big(\psi(x,\,\cdot\,)\big)_{x\in
  S_X}$ are distribution determining. In this case, the transition
probabilities of $X$ and $Y$ determine each other uniquely through
(\ref{Markdual}).

In our case, each of the voter models $X$ we have introduced is a cancellative
spin system in the sense of \cite{Gri79}, hence it is know that for each
of these models, there exists a Markov process $Y$ which is also a cancellative
spin system and which is dual to $X$ with the duality function
$\psi(x,y):=(-1)^{|xy|}$, where for any $x,y,z\in\{0,1\}^\Z$ we let
$xy(i):=x(i)y(i)$ denote the componentwise product of $x$ and $y$, and we
write $|z|:=\sum_iz(i)$. In this case, (\ref{Markdual}) implies that
\be\label{dual}
\P\big[|X_tY_0|\mbox{ is odd}\big]=\P\big[|X_0Y_t|\mbox{ is odd}\big]
\qquad(t\geq 0)
\ee
whenever $X$ and $Y$ are independent (with arbitrary initial laws), where in
order for our expressions to be well-defined we must assume that either
$|X_0|$ or $|Y_0|$ is finite.

A special property of the rebellious voter model (that motivated its
introduction in \cite{SS08hv}) is that its interface model and dual model
coincide, i.e., if $X$ is the rebellious voter model, then (\ref{dual}) holds
for the same $Y$, with dynamics described in (\ref{twosidY}), that also
describes the interfaces of $X$. One way of checking this is to verify
(\ref{Gdual}); an alternative proof, based on the graphical representation of
cancellative spin systems, is described at \cite[formula~(1.8)]{SS08hv}. For the
one-sided rebellious voter model, the dual and interface models do not
coincide, but its dual is the mirror image of its interface model, i.e., the
one-sided rebellious voter model satisfies (\ref{dual}) for a model $Y$ whose
dynamics are described by (compare (\ref{onesidY}))
\be\label{mirrorY}
y\mapsto y^{\{i-1,i\}}\quad\mbox{with rate}
\quad\al1_{\{y(i)=1\}}+(1-\al)1_{\{y(i+1)=1\}}.
\ee
It turns out that the dual (in the sense of cancellative spin systems) of the
disagreement model is the DBARW (which is the interface model of the swapping
voter model) while the dual of the swapping voter model is the SARW (which is
the interface model of the disagreement voter model) (see
\cite[Section~2.1]{SS08hv}). Note that in the SARW there is no branching,
which has far-reaching consequences for both the disagreement and swapping
voter models.

\subsection{Coexistence, survival and interface tightness}

It is easy to see that the constant configurations $\un 0$ and $\un 1$ are
traps for the various voter models we have introduced. By definition, we say
that a Markov process taking values in $\{0,1\}^{\Z^d}$ exhibits {\em
  coexistence} if there exists an invariant law that is concentrated on
configurations that are not identically zero or one. It is known that
nearest-neighbor voter models on $\Z^d$ exhibit coexistence if and only if
$d>2$; in particular, for $\al=1$, our models, being one-dimensional voter
models, do not exhibit coexistence. On the other hand, our models exhibit
coexistence for $\al=0$ since in this case one can check that product measure
with intensity a half is an invariant law.

The question is whether there is a nontrivial phase transition between
coexistence and noncoexistence. For the disagreement and swapping voter models
it turns out that this is not the case. Indeed, it can be rigorously proved
\cite{Sud90,NP99,SS08hv} that these model exhibits noncoexistence for all
$\al>0$. On the other hand, it appears that the rebellious voter model and its
one-sided counterpart do exhibit a nontrivial phase transition between
coexistence and noncoexistence at some $0<\al_{\rm c}<1$.

We will mainly be interested in two functions of our processes. It is said
that a voter model {\em survives} if there is a positive probability that the
process started with a single one never gets trapped in $\un 0$. In view of
this, for the (one- or two-sided) rebellious voter model with competition
parameter $\al$, we will be interested in the {\em survival probability}
\be\label{rhodef}
\rho(\al):=\P^{\de_0}\big[X_t\neq\un 0\ \forall t\geq 0\big],
\ee
where $\P^{\de_0}$ denotes the law of the process started in $X_0=\de_0$,
with $\de_0(i):=1_{\{i=0\}}$. It seems plausible that for many models of the
type we are considering, survival is equivalent to coexistence. Indeed, using
the fact that its interface and dual models coincide, this can be rigorously
proved for the rebellious voter model \cite[Lemma~2]{SS08hv}; the
same proof works for the one-sided process. One of the main aims of the
present paper is to obtain good numerical data for the function $\rho$.

The function $\rho$ quantifies, in a sense, how strongly the process exhibits
coexistence (see also formulas (\ref{inter}) and (\ref{approx}) below). We
need a similar formula to quantify noncoexistence. For any $i\in\Z$, let
$x^{\rm H}_i\in\{0,1\}^\Z$, defined as
\be
x^{\rm H}_i(j):=\left\{\ba{ll}
0\quad&\mbox{if }j<i,\\
1\quad&\mbox{if }i\leq j
\ea\right.
\ee
denote a `Heaviside' configuration of zeros and ones that are `completely
separated'. Following \cite{CD95}, we say that a model $X$ (with given parameter
$\al$) exhibits {\em interface tightness} if
\be\label{chidef}
\chi(\al):=\lim_{t\to\infty}\P^{x^{\rm H}_0}
\big[X_t=x^{\rm H}_i\mbox{ for some }i\in\Z\big]>0,
\ee
i.e., if the process started in a Heaviside configuration spends a positive
fraction of its time in such configurations. It is intuitively plausible that
interface tightness implies noncoexistence. It is rigorously known that the
swapping voter model exhibits interface tightness for all $\al>0$
\cite{SS08if}.

The main aim of the present paper is to obtain good numerical data for the
functions $\rho$ and $\chi$, both for the rebellious voter model and for its
one-sided analogue.

\section{Main results}

\subsection{Methods}

In our simulations, we start the interface process $Y$ of either the
rebellious voter model or its one-sided counterpart with an odd number of ones
on an interval of $N$ sites with periodic boundary conditions. Note that in
this case, because the number of ones is odd and because of the periodic
boundary conditions, we can no longer represent $Y$ as in (\ref{inter}), but
the dynamics in (\ref{onesidY}) and (\ref{twosidY}), as well as the duality
relation (\ref{dual}), still make sense. Since we start with one particle,
because of parity conservation, the system cannot die out. Letting $Y_\infty$
denote the process in equilibrium, we assume that for large $N$, the following
approximations are valid:
\be\ba{rl}\label{approx}
{\rm(i)}&\rho(\al)\approx 2\P\big[Y_\infty(0)=1\big]
=\ffrac{2}{N}\E\big[|Y_\infty|\big],\\[5pt]
{\rm(ii)}&\chi(\al)\approx \P\big[|Y_\infty|=1\big].
\ec
The arguments why formula (\ref{approx})~(ii) should hold are obvious from
(\ref{inter}) and (\ref{chidef}). The justification of formula
(\ref{approx})~(i) is more subtle. Let $Y^{1/2}$ denote the interface process
of the rebellious voter model on $\Z$, started in an initial law such that the
$(Y_0(i))_{i\in\Z}$ are independent with
$\P[Y_0(i)=0]=\P[Y_0(i)=1]=\frac{1}{2}$, and let $X$ be the rebellious voter
model started in $X_0=\de_0$. Then by duality (\ref{dual}),
\be\label{dusur}
\P\big[Y^{1/2}_t(0)=1\big]=\P\big[|Y^{1/2}_tX_0|\mbox{ is odd}\big]
=\P\big[|Y^{1/2}_0X_t|\mbox{ is odd}\big]=\ffrac{1}{2}\P\big[X_t\neq\un{0}\big],
\ee
where $\un{0}$ denotes the constant zero configuration. It is known that the
law of $Y^{1/2}_t$ converges as $t\to\infty$ to an invariant law, which is in
fact the unique spatially homogeneous, nontrivial invariant of the process
\cite[formula (1.12) and Theorem~3~(b)]{SS08hv}. Taking the limit $t\to\infty$
in (\ref{dusur}), we see that for this invariant law, formula
(\ref{approx})~(i) holds with equality. Since it seems reasonable that the
invariant laws of our finite systems approximate the invariant law of the
infinite system, this justifies (\ref{approx})~(i). Similar arguments apply to
the one-sided model.

\subsection{A first qualitative look}

\begin{figure}
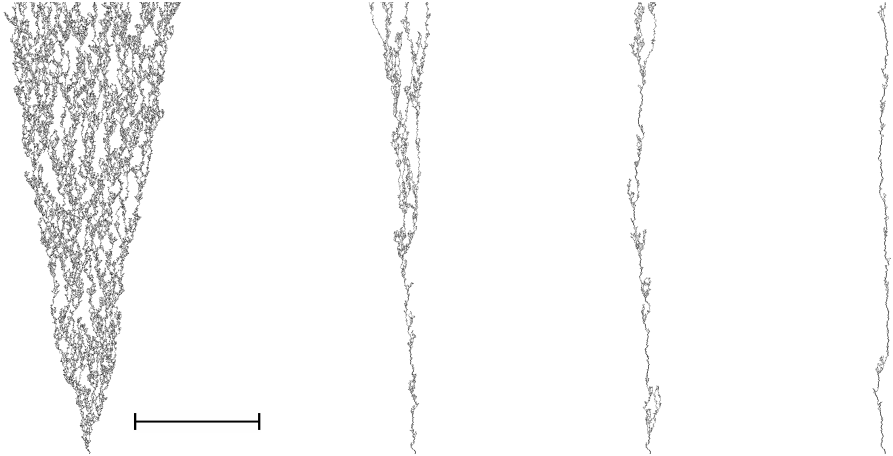

\begin{center}
\includegraphics[width=4cm,height=6cm]{\figfilebw{vier500a}}
\hspace{.3cm}
\includegraphics[width=3cm,height=6cm]{\figfilebw{vijf}}
\includegraphics[width=3cm,height=6cm]{\figfilebw{criti}}
\includegraphics[width=3cm,height=6cm]{\figfilebw{zes}}
\caption[The interface model of the rebellious voter model for $\al=0.4$,
  $\al=0.5$, $\al=0.51$ and $\al=0.6$.] {The interface model $Y$ of the
  rebellious voter model for $\al=0.4$, $\al=0.5$, $\al=0.51$ and $\al=0.6$,
  started with a single interface. Space is plotted horizontally and time
  vertically, with a black spot at a point $(i,t)$ indicating that
  $Y_t(i)=1$. Total time elapsed in each picture is 1800. The horizontal bar is 500 sites long.}
\label{fig:bitmap}
\end{center}
\end{figure}

Before we present our numerical data for the functions $\rho$ and $\chi$ from
(\ref{approx}), we first take a look at the qualitative behavior of the
interface model $Y$ of the rebellious voter model and its one-sided analogue.
In Figure~\ref{fig:bitmap}, we have graphically plotted $Y$ for the two-sided
process, started with a single one, with time running upwards and a black spot
at a point $(i,t)$ indicating that $Y_t(i)=1$. The figures shows the process
for four values of $\al$, where $\al_{\rm c}\approx 0.51$ is the estimated
critial value for this process (see Section~\ref{S:twosid} below).

The pictures show that for $\al<\al_{\rm c}$, all ones are contained in an
interval with edges that grow at an approximately linear speed, and that
inside this interval, ones occur at some approximately constant equilibrium
density. On the other hand, for $\al>\al_{\rm c}$, the process spends a
positive fraction of its time in states where there is just a single site in
state one, which indicates that the rebellious voter model exhibits interface
tightness (see (\ref{chidef})). What happens exactly at $\al=\al_{\rm c}$ is
less clear from these pictures, but our numerical data, presented in
Sections~\ref{S:twosid} and \ref{S:edge} below, suggest that at $\al=\al_{\rm
  c}$, the edge speed is zero and the rebellious voter model exhibits neither
coexistence nor interface tightness.

\begin{figure}
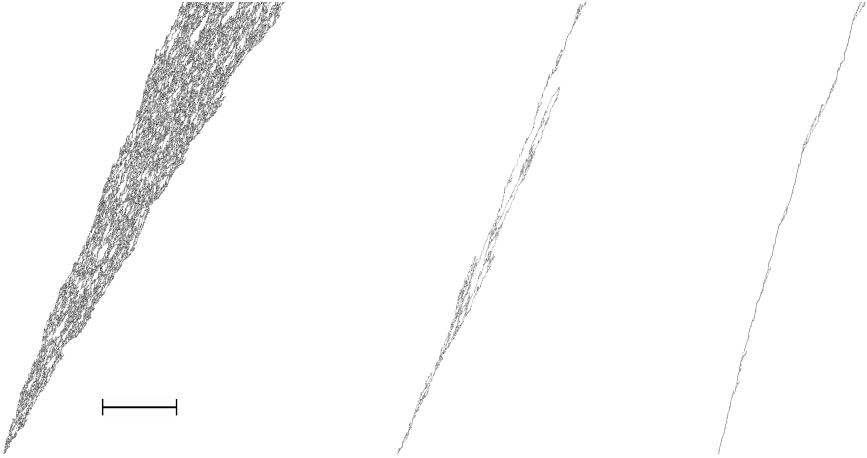

\begin{center}
\includegraphics[width=4cm,height=6cm]{\figfilebw{mdrie500}}
\hspace{1cm}
\includegraphics[width=3cm,height=6cm]{\figfilebw{mvijf}}
\hspace{1cm}
\includegraphics[width=3cm,height=6cm]{\figfilebw{mzes}}
\caption[The interface model of the one-sided rebellious voter model for
  $\al=0.3$, $\al=0.5$ and $\al=0.6$.]{The interface model of the one-sided
  rebellious voter model for $\al=0.3$, $\al=0.5$ and $\al=0.6$. Total time elapsed in each picture is 1800. The horizontal bar is 500 sites long.}
\label{fig:bitmap2}
\end{center}
\end{figure}

The picture for the one-sided process is very similar, except that now, due to
the nature of the process, the process as a whole tends to drift to the right
at some approximately constant speed; see Figure~\ref{fig:bitmap2}.

\subsection{The two-sided model}\label{S:twosid}

We simulated the interface model $Y$ of the rebellious voter model on an
interval of $N$ sites with periodic boundary conditions by slowly increasing
or decreasing $\al$ from some initial value $\al_{\rm b}$ to some final value
$\al_{\rm e}$ during a time interval of length $T$. We then divided our time
interval into $n$ equal pieces and plotted the average value of $2|Y|/N$ (for
the function $\rho$) or the fraction of the time that $|Y|=1$ (for the
function $\chi$) against the average value of $\al$, for each of the $n$ time
intervals.

This method has both advantages and disadvantages. One of its main advantages
is that it allows one to quickly obtain data for a very large number of values
of the parameter, which is useful if one is interested in curve fitting or in
estimating derivative functions such as in Figure~\ref{fig:difharm}. The
method introduces obvious errors due to the nonzero speed of varying the
parameter, but it is not difficult to get a rough idea of the size of these
effects. In fact, one can get a good idea of the size of the errors just by
looking at the curves (see Section~\ref{S:finsiz} for a more detailed
discussion). A disadvantage of the method is that the size of the errors
varies a lot along the curves since relaxation times increase as one
approaches the critical point. One can partially compensate for this by doing
detailed simulations near the critical point. (A better approach, which we
have not pursued, would be to vary $\al$ at a suitable, precisely chosen
nonconstant speed.)

The statistical errors in our data for the function $\rho$ appear to be
approximately normally distributed (with a variance that depends strongly on
$\al$), but those in the data for $\chi$ have a skew distribution that is very
far from a normal one (this is somewhat visible from
Figure~\ref{fig:oneserr}). This non-Gaussian behavior is not
unexpected. Indeed, if we start the interface model $Y$ in a state with
$|Y_0|=3$ particles, then it is known that for the model with $\al=1$,
\be
\P\big[\inf\{s\geq 0:|Y_s|=1\}\geq t\big]\sim t^{-3/2}\quad\mbox{as }t\to\infty.
\ee
It is seems natural to conjecture that this is also true for $\al<1$, hence
the process started with a single particle makes excursions into states with
more than one particle that have a duration which has a finite first moment
but infinite second moment.

\begin{figure}
\begin{center}
\includegraphics[width=12cm,height=6cm]{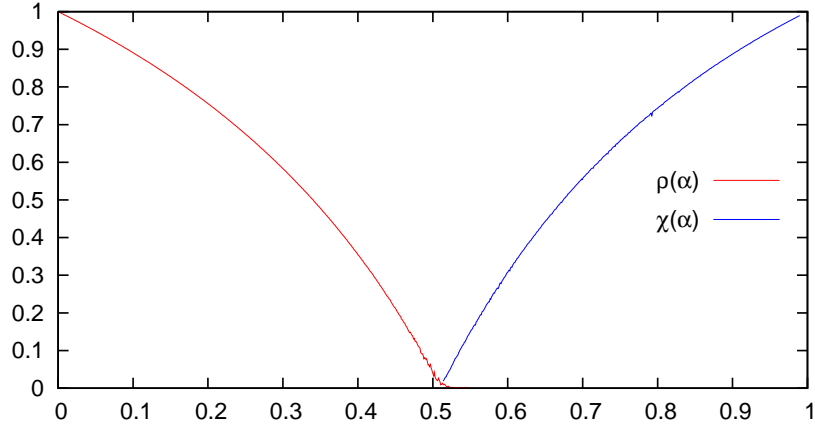}
\caption[The functions $\rho$ and $\chi$ for the rebellious voter model.]
{The functions $\rho$ and $\chi$ for the rebellious voter model. Data obtained
  with parameters as in (\ref{sympar}).}
\label{fig:double}
\end{center}
\end{figure}

In Figure~\ref{fig:double} we have plotted approximations for the
functions $\rho$ and $\chi$ obtained by our methods, using the parameters:
\be\ba{l@{\quad}l}\label{sympar}
\mbox{for }\rho:&N=2^{15},\ n=2^9,\ T=10^8,
\ \al_{\rm b}=0,\ \al_{\rm e}=0.55,\\
\mbox{for }\chi:&N=2^{15},\ n=2^9,\ T=10^{12},
\ \al_{\rm b}=0.99,\ \al_{\rm e}=0.49.
\ec

The values $\rho(0)=1$ and $\chi(1)=1$ are expected on theoretical grounds.
Indeed, it is easy to check that the rebellious voter model with $\al=0$ never
dies out (recall (\ref{rhodef})) while in the pure voter case $\al=1$, the
process started in a Heaviside state stays in such a state for all time
(recall (\ref{chidef})).

\begin{figure}
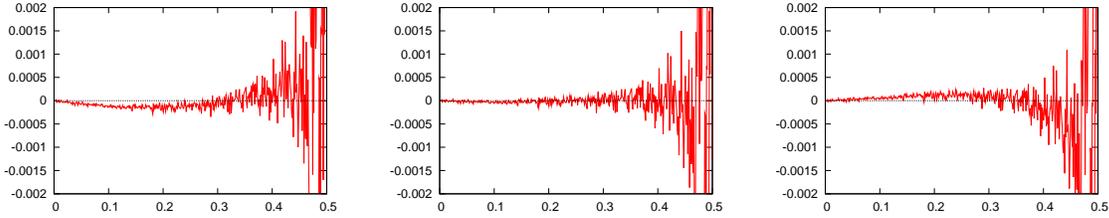

\begin{center}
\includegraphics[width=5cm]{\figfile{twoerr_min2}}
\includegraphics[width=5cm]{\figfile{twoerr2}}
\includegraphics[width=5cm]{\figfile{twoerr_plus2}}
\caption[Comparison of $\rho$ with different linear fractional functions.]
{The function $\rho(\al)-(1-c_1\al)/(1-c_2\al)$ for the values $c_1=1.957$,
  $c_2=0.973$ (left), $c_1=1.958$, $c_2=0.975$ (middle) and $c_1=1.959$,
  $c_2=0.977$ (right).}
\label{fig:twoerr}
\end{center}
\end{figure}

\begin{figure}
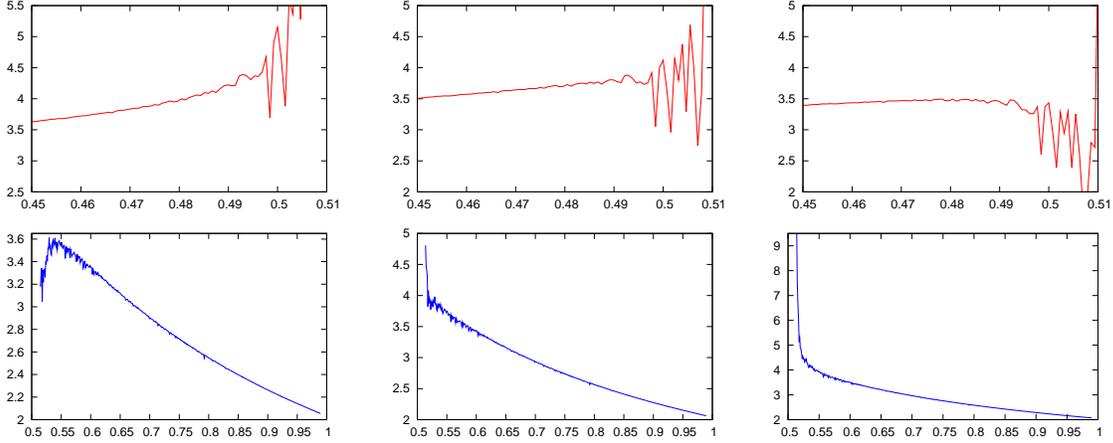

\begin{center}
\includegraphics[width=5cm]{\figfile{cridet5}}
\includegraphics[width=5cm]{\figfile{cridet1}}
\includegraphics[width=5cm]{\figfile{cridet2}}
\includegraphics[width=5cm]{\figfile{upturn508}}
\includegraphics[width=5cm]{\figfile{upturn51}}
\includegraphics[width=5cm]{\figfile{upturn512}}
\caption[Estimation of the critical point.]  {Estimation of the critical point
  of the rebellious voter model. Above: The function $\rho(\al)/(\al_0-\al)$
  for the values $\al_0=0.508$ (left), $\al_0=0.510$ (middle) and
  $\al_0=0.512$ (right). Below: the function $\chi(\al)/(\al-\al_0)$ for the
  values $\al_0=0.508$ (left), $\al_0=0.510$ (middle) and $\al_0=0.512$
  (right). Our simulations for $\rho$ use the parameters $N=2^{15}$, $n=2^7$,
  $T=10^9$, $\al_{\rm b}=0.55$, $\al_{\rm e}=0.45$ and
  those for $\chi$ use the parameters in (\ref{sympar}).}
\label{fig:upturn}
\end{center}
\end{figure}

The numerical data for $\rho$ are well fitted by a linear fractional function
of the form
\be\label{rhofrac}
\rho(\al)=\frac{1-c_1\al}{1-c_2\al}\quad\mbox{with}\quad
c_1=1.958\pm 0.001\mbox{ and }c_2=0.975\pm 0.002,
\ee 
(see Figure~\ref{fig:twoerr}). Assuming this sort of function fitting is
justified, we arrive at an estimate for the critical point of
\be
\al_{\rm c}=0.5107\pm 0.0003.
\ee
This estimate is to be viewed with caution, since this is based on
extrapolation of data for (approximately) $0\leq\al\leq 0.35$, assuming that
the linear fractional form (\ref{rhofrac}) holds for all $\al\leq\al_{\rm
  c}$. A more robust method (see Figure~\ref{fig:upturn}) based on our best
data for the functions $\rho$ and $\chi$ near the critical point leads to the
estimate:
\be
\al_{\rm c}=0.510\pm 0.002.
\ee
(See also Figure~\ref{fig:detail}.) It seems that linear fractional functions
do not fit the numerical data for $\chi$ really well.

Our data strongly suggest that $\rho(\al_{\rm c})=0=\chi(\al_{\rm c})$, which
implies that at criticality, the process exhibits neither coexistence nor
interface stability. Summarizing, the picture that emerges from our numerical
data is as follows:
\begin{quote}
There exists a critical value $\al_{\rm c}\approx 0.510$ such that the
rebellious voter model exhibits coexistence if and only if $\al<\al_{\rm c}$
and interface tightness if and only if $\al>\al_{\rm c}$. The function $\rho$
from (\ref{rhodef}) has a finite negative slope and is concave on $[0,\al_{\rm
    c}]$, and is approximately given by a linear fractional function of the
form (\ref{rhofrac}). The function $\chi$ from (\ref{chidef}) has a finite
positive slope and is concave on $[\al_{\rm c},0]$.
\end{quote}

\subsection{The one-sided model}

\begin{figure}
\begin{center}
\includegraphics[width=12cm,height=6cm]{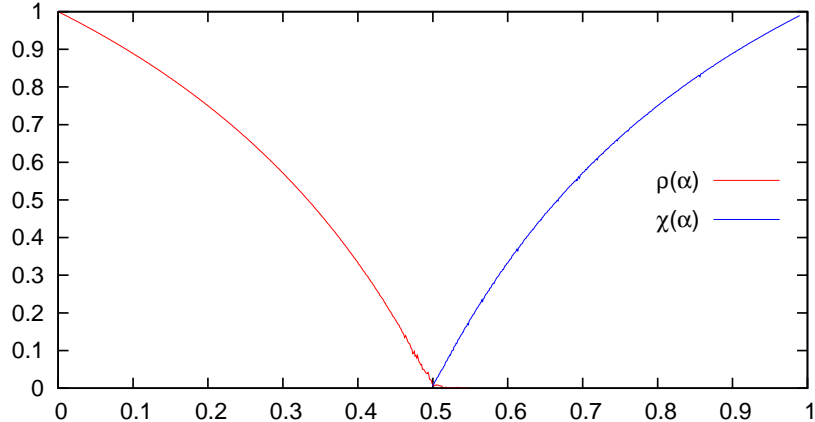}
\caption[The functions $\rho$ and $\chi$ for the one-sided rebellious voter
  model.]{The functions $\rho$ and $\chi$ for the one-sided
  rebellious voter model. Data obtained with parameters as in (\ref{sympar}).}
\label{fig:onesid}
\end{center}
\end{figure}

We have run the same simulation as described in the previous section, with
parameters as in (\ref{sympar}), also for the one-sided rebellious voter
model. The resulting approximations for the functions $\rho$ and $\chi$ are
plotted in Figure~\ref{fig:onesid}.

\begin{figure}
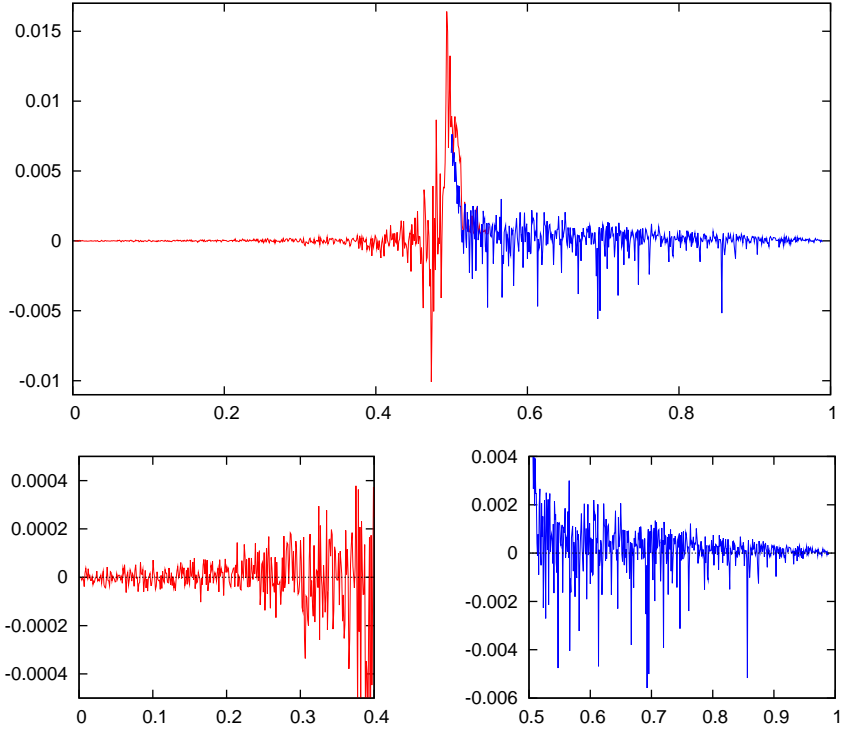

\begin{center}
\includegraphics[width=12cm,height=6cm]{\figfile{oneserr2}}
\includegraphics[width=6cm,height=4cm]{\figfile{Loneserr2}}
\includegraphics[width=6cm,height=4cm]{\figfile{Roneserr2}}
\caption[Difference between explicit formulas and the functions $\rho$ and
  $\chi$ for the one-sided rebellious voter model.]{Difference between the
  explicit formulas in (\ref{explicit}) and the functions $\rho$ and
  $\chi$ for the one-sided rebellious voter model, with detail of the
  left and right side of the picture.}
\label{fig:oneserr}
\end{center}
\end{figure}

\begin{figure}
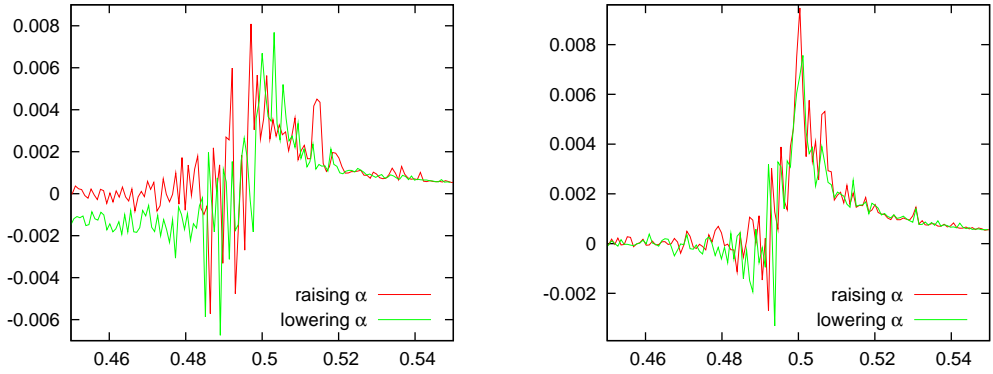

\begin{center}
\includegraphics[width=7cm]{\figfile{compar_fast}}
\includegraphics[width=7cm]{\figfile{compar_slow}}
\caption[Difference between explicit formula and the functions $\rho$ for the
  one-sided rebellious voter model.]{Plots of
  $\rho(\al)-\max\{0,(1-2\al)/(1-\al)\}$ for the one-sided rebellious voter
  model near the critical point. All plots use the parameters $N=2^{15}$,
  $n=2^7$. On the left: $T=10^9$ and $\al_{\rm b}=0.445$,
  $\al_{\rm e}=0.55$ resp.\ $\al_{\rm b}=0.55$ and $\al_{\rm e}=0.45$. 
  On the right: $T=4\cdot 10^9$ and $\al_{\rm b}=0.445$,
  $\al_{\rm e}=0.55$ resp.\ $\al_{\rm b}=0.55$ and $\al_{\rm e}=0.445$.}
\label{fig:compar}
\end{center}
\end{figure}

At first sight, Figure~\ref{fig:onesid} looks extremely similar to
Figure~\ref{fig:double}. A closer inspection reveals, however, that the two
plots are not identical. Indeed, it seems that for the one-sided model, the
functions $\rho$ and $\chi$ are described by the explicit formulas:
\be\label{explicit}
\rho(\al)=\dis 0\vee\frac{1-2\al}{1-\al}
\quad\mbox{and}\quad
\chi(\al)=0\vee\big(2-\frac{1}{\al}\big).
\ee
In particular, one has the symmetry $\rho(1-\al)=\chi(\al)$ and the critical
parameter seems to be given by
\be
\al_{\rm c}=0.500\pm 0.002\quad =\ffrac{1}{2}\quad(?)
\ee
In Figure~\ref{fig:oneserr}, we have plotted the differences between the
explicit functions in (\ref{explicit}) and our data for $\rho$ and $\chi$,
respectively. The systematic deviations of $\rho$ and $\chi$ from the proposed
formulas near the critical point probably stem from finite size effects
(compare Figures~\ref{fig:finN} and \ref{fig:finchi} below). A more detailed
comparison of $\rho$ with its proposed explicit formula is shown in
Figure~\ref{fig:compar}. Because of the small difference between the left and
right edge speed near the critical point (see Section~\ref{S:edge}), these
detailed simulations are sensitive to the direction in which $\al$ is varied
(compare Figure~\ref{fig:hyster} and the discussion in
Section~\ref{S:finsiz}). In particular, the systematic deviation between our
proposed explicit formula and the curve produced by lowering $\al$ in the left
plot of Figure~\ref{fig:compar} seems to stem from such effects. This
systematic deviation almost disappears after giving the system more time to
relax (see the right plot of Figure~\ref{fig:compar}).

\begin{figure}
\begin{center}
\includegraphics[width=9cm]{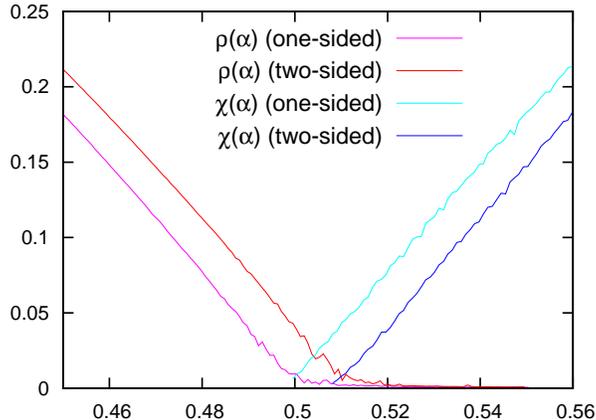}
\caption[Detail near the critical point.]  {A detailed simulation near the
  critical point of the functions $\rho$ and $\chi$ for the two-sided process
  and the one-sided process. On the left: our simulations
  for $\rho$ use the parameters $N=2^{15}$, $n=2^7$, $T=4\cdot10^9$, $\al_{\rm
    b}=0.445$, $\al_{\rm e}=0.55$. Our plots for $\chi$ show detail of a
  simulation with the parameters $N=2^{15}$, $n=2^9$, $T=10^{12}$, $\al_{\rm
    b}=0.99$, $\al_{\rm e}=0.5$.}
\label{fig:detail}
\end{center}
\end{figure}

We guessed our formula for $\rho$ after estimating the first few derivatives
at $\al=0$ and the formula for $\chi$ by analogy with $\rho$. We do not know
of any theoretical reason to expect these formulas. A detailed
comparison of the two-sided and one-sided model near the critical point is
shown in Figure~\ref{fig:detail}.

\begin{figure}
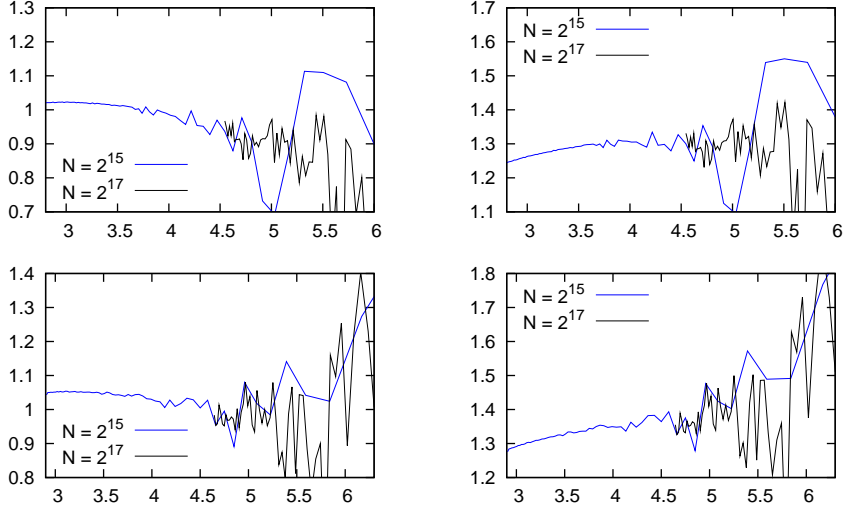

\begin{center}
\includegraphics[width=6cm, height=3.5cm]{\figfile{elog92}}
\includegraphics[width=6cm, height=3.5cm]{\figfile{elog100}}
\includegraphics[width=6cm, height=3.5cm]{\figfile{emlog92}}
\includegraphics[width=6cm, height=3.5cm]{\figfile{emlog100}}
\caption[Double logarithmic plots.]{Estimation of the order parameter critical
  exponent $\bet$. Plots of $\log(\rho)-\bet\log(\al_{\rm c}-\al)$ as a
  function of $-\log(\al_{\rm c}-\al)$. For the right choice of $\bet$, the
  asymptotic slope of these functions should be zero as the parameter on the
  horizontal axis tends to infinity. Above: plots for the two-sided model with
  $\bet=0.92$ (left) and $\bet=1$ (right), based on the estimate $\al_{\rm
    c}=0.51$. Below: the same for the one-sided model, based on the estimate
  $\al_{\rm c}=0.5$. Each figure combines two plots, one based on simulations
  with the parameters $N=2^{15}$, $n=2^7$, $T=4\cdot10^9$, $\al_{\rm
    b}=0.445$, $\al_{\rm e}=0.55$, the other instead using $N=2^{17}$ and a
  smaller range of alpha ($\al_{\rm b}=0.5$, $\al_{\rm e}=0.52$ for the
  two-sided model and $\al_{\rm b}=0.49$, $\al_{\rm e}=0.51$ for the one-sided
  model).}\label{fig:logplot}
\end{center}
\end{figure}

Our proposed explicit formula for $\rho$, if it is correct, implies that the
order parameter critical exponent $\bet$ for our model is one. Since this is
different from the recently published values $\bet=0.92(2)$ \cite{Hin00} and
$\bet=0.95(1)$ \cite{OM06} for the PC universality class, we have used our
data to obtain a direct estimate for the value of $\bet$. A plot of
$\log(\rho)$ as a function of $\log(\al_{\rm c}-\al)$ yields an approximate
linear graph as $\al$ approximates $\al_{\rm c}$ from below, the slope of
which should be $\bet$. To estimate this slope, in Figure~\ref{fig:logplot} we
have plotted $\log(\rho)-\bet\log(\al_{\rm c}-\al)$ as a function of
$-\log(\al_{\rm c}-\al)$, both for the two-sided and one-sided model, for
different values of $\bet$. These plots are suggestive of a value of $\bet$
that is perhaps closer to $0.92$ than to $1$, but this is probably due to the
fact that our simulations are not precise enough to get sufficiently close to
the critical point. In fact, our data become unreliable at $\al_{\rm
  c}-\al\leq 0.004$, which corresponds to the point $-\log(0.004)\approx 5.5$ in
Figure~\ref{fig:logplot}. Note that in Figure~\ref{fig:compar}, for the
one-sided model, the same data were shown to be consistent with our
proposed explicit formula for $\rho$ in (\ref{explicit}), which, if true,
implies that $\bet=1$. We conclude that our data are not good enough to
convincingly distinguish between $\bet=0.92$ and $\bet=1$.

\subsection{Finite size effects}\label{S:finsiz}

\begin{figure}
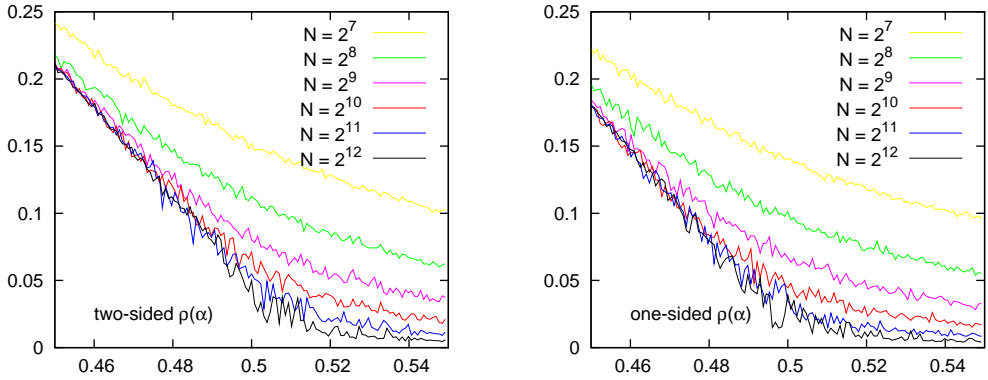

\begin{center}
\includegraphics[width=7cm]{\figfile{symfinN}}
\includegraphics[width=7cm]{\figfile{finN2}}
\caption[Effect of the system size on our approximation of $\rho$.]  {Effect
  of the system size on our approximation of the function $\rho$ from
  (\ref{approx}). Plotted are our approximations for $\rho$ using the
  parameters $n=2^7$, $T=10^8$, $\al_{\rm b}=0.55$, $\al_{\rm e}=0.45$, and
  the system sizes $N=2^7,2^8,2^9,2^{10},2^{11},2^{12}$, respectively. On the
  left: the two-sided model. On the right: the one-sided model.}
\label{fig:finN}
\end{center}
\end{figure}

To gain more insight into how close our numerical data are to the real
functions $\rho$ and $\chi$, we tested the effect of varying the system size
$N$. In Figure~\ref{fig:finN} we have plotted our approximation (see
(\ref{approx})) of the function $\rho$ for different values of the system size
$N$. The pictures for the two-sided and one-sided process are extremely
similar, except for a shift of the critical point. In view of
Figure~\ref{fig:finN}, it seems that we can rule out the possibility that the
observed differences between the two-sided and one-sided process are due to
finite space effects (and a hypothetical slower convergence for the
two-sided model).

Note that in Figure~\ref{fig:finN}, our approximations to $\rho$ near the
critical point become more rough as the system size is increased. This can be
understood due to two effects. On the one hand, near criticality, the random
variable $|Y_\infty|$ assumes the values $1,3,5,7,\ldots$ with approximately
equal probabilities (see Section~\ref{S:freq}), which means that the number of
ones in the system fluctuates from being close to one to a positive fraction
of all sites in the lattice. As the system size gets large, this means huge
fluctuations, which are moreover slow since at criticality, the edge speed is
zero.

\begin{figure}
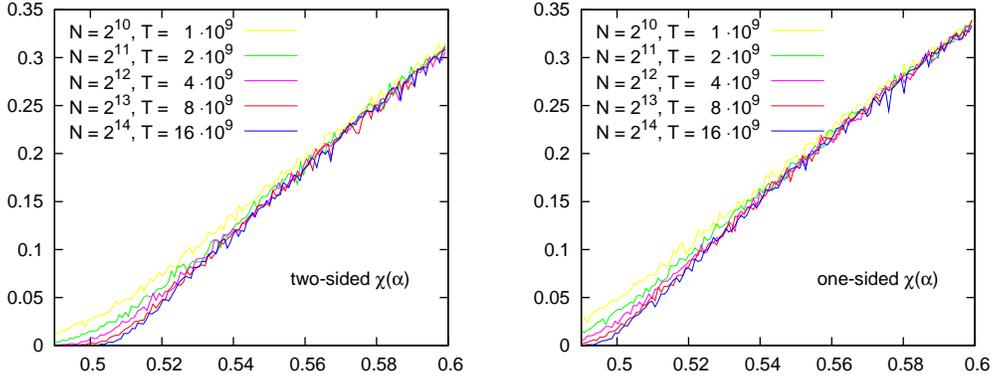

\begin{center}
\includegraphics[width=7cm]{\figfile{finN3}}
\includegraphics[width=7cm]{\figfile{finN4}}
\caption[Effect of the system size on our approximation of $\chi$.]  {Effect
  of the system size on our approximation of the function $\chi$ from
  (\ref{approx}). Plotted are our approximations for $\chi$ using the
  parameters $n=2^7$, $\al_{\rm b}=0.6$, $\al_{\rm e}=0.49$--$0.5$, and the
  system sizes $N=2^{10},2^{11},2^{12},2^{13},2^{14}$ and times
  $T=10^9,2\cdot10^9,4\cdot10^9,8\cdot10^9,16\cdot10^9$, respectively. On the
  left: the two-sided model. On the right: the one-sided model.}
\label{fig:finchi}
\end{center}
\end{figure}

In Figure~\ref{fig:finchi}, we have plotted our approximations for the
function $\chi$ for various values of $N$. To counteract the `roughening'
effect we have just described, in these simulations, we have increased the
total time $T$ together with the system size $N$. Again, the pictures for
the two-sided and one-sided process are similar except for a small shift in
the critical point.

\begin{figure}
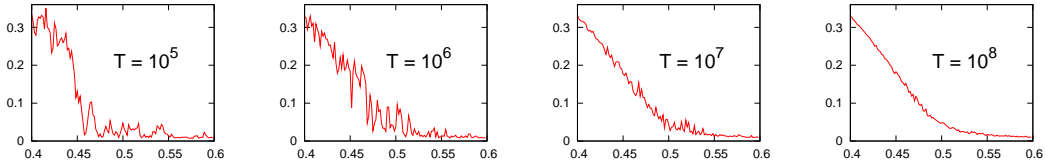

\begin{center}
\includegraphics[width=3.5cm]{\figfile{tijd5}}
\includegraphics[width=3.5cm]{\figfile{tijd6}}
\includegraphics[width=3.5cm]{\figfile{tijd7}}
\includegraphics[width=3.5cm]{\figfile{tijd8}}
\caption[Effect of the simulation time on our approximation of $\rho$.]
        {Effect of the duration of our simulations on our approximation of the
          function $\rho$ from (\ref{approx}) for the one-sided model. Plotted
          are our approximations for $\rho$ using the parameters $N=2^{10}$,
          $n=2^7$, $\al_{\rm b}=0.6$, $\al_{\rm e}=0.4$, and the total elapsed
          times $T=10^5,10^6,10^7,10^8$, respectively.}
\label{fig:finT}
\end{center}
\end{figure}

\begin{figure}
\begin{center}
\includegraphics[width=8cm]{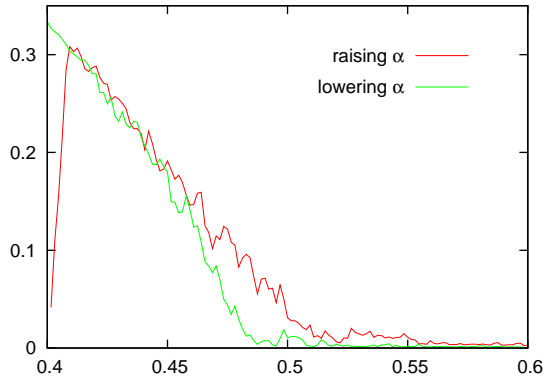}
\caption[Effects of choosing the time too small compared to the system size.]
        {Effects of choosing the time too small compared to the system size,
          in simulations of $\rho$ for the one-sided model. The two curves
          were produced by raising resp.\ lowering $\alpha$ between the values
          $\al_{\rm b}$ resp.\ $\al_{\rm e}=0.4$ and $\al_{\rm e}$
          resp.\ $\al_{\rm b}=0.6$. The other simulation parameters were
          $n=2^7$, $N=2^{13}$, and $T=10^6$.}
\label{fig:hyster}
\end{center}
\end{figure}

We have also run simulations for different values of the total elapsed time,
to investigate the effect of this parameter on the quality of our data.  Since
short times mean the system does not have enough time for time averages to
reach their equilibrium values, the main effect of increasing $T$ is to
smoothen our approximated functions. For the function $\rho$, this effect is
demonstrated in Figure~\ref{fig:finT}. A second effect of choosing a too short
time $T$ is that the numerical functions `lag behind' in that they show values
belonging to an $\al$ that lies somewhat in the past. This effect is
demonstrated in Figure~\ref{fig:hyster}. In this case, our approximations near
the critical point depend on the direction in which we vary $\al$: the graph
produced by increasing $\al$ `overshoots' the critical point while the graph
produced by lowering $\al$ picks up too late. Also note the `hook' at the
beginning of the graph started with $\al$ below the critical point, which is
due to the fact that we start with a single one and the finite edge speed
needs several time steps to fill out all space. This sort of effects only
occur when $T$ is small relative to $N$; in all of our simulations (except
those in Figure~\ref{fig:hyster}) we have done our best to choose our
parameters such that these effects in minimal. In case of doubt, we have run
the same simulation in both directions for comparison (see
Figure~\ref{fig:compar}). In some plots, we have removed a small `hook' at the
beginning.

\section{Other functions of the process}

\subsection{Harmonic functions}\label{S:harm}

We have obtained numerical data for a number of other functions of our
processes. We have mainly concentrated on the one-sided model, in view of the
apparent explicit formulas for the functions $\rho$ and $\chi$ in this
case. We have not been able, however, to find explicit formulas for any other
functions than these two. Nevertheless, our additional functions show some
interesting phenomena, which we discuss here. In the present section, we
discuss a harmonic function for the (one-sided) rebellious voter model. In the
next two sections, we look at the equilibrium probabilities of seeing
$3,5,7,\ldots$ ones in the interface model, and edge speeds, respectively.

Recall that if $X$ is a Markov process with state space $S_X$, then a function
$f:S_X\to\R$ is called a {\em harmonic function} for $X$ if, for any initial
state, the process $(f(X_t))_{t\geq 0}$ is a martingale. At least for finite
state spaces, this is equivalent to the statement that $G_Xf=0$, where $G_X$
denotes the generator of $X$. In general, a duality between Markov processes
translates the invariant laws of one process into harmonic functions of the
other process. More precisely, if $X$ and $Y$ are dual Markov process,
with duality function $\psi$, and $\mu$ is an invariant measure for $Y$, then
using (\ref{dual}) it is easy to see that setting
\be
f(x):=\int_{S_Y}\psi(x,y)\mu(\di y)\qquad(x\in S_X)
\ee
defines a harmonic function $f$ for $X$. Note that multiplying $f$ with a
constant will not change the fact that it is harmonic.

To apply this to the one-sided rebellious voter model, we first need to
introduce some notation. In this section, we let $X_t$ denote the rebellious
voter model, $\rvec X_t$ the one-sided rebellious voter model, and $\lvec X_t$
the mirror image of the latter, i.e., the model that jumps as
\be\ba{l}\label{mirRV}
x\mapsto x^{\{i\}}\quad\mbox{with rate}\\[5pt]
\dis\quad\al 1_{\{x(i)\neq x(i+1)\}}+(1-\al)1_{\{x(i+1)\neq x(i+2)\}}.
\ec
We let $\rvec Y_t$ and $\lvec Y_t$ denote the interface models of $\rvec X_t$
and $\lvec X_t$, respectively, i.e., $\rvec Y_t$ is the process with dynamics
described in (\ref{onesidY}) and $\lvec Y_t$ is the process in
(\ref{mirrorY}). Then $\lvec Y_t$ is dual to $\rvec X_t$ and $\rvec Y_t$ is
dual to $\lvec X_t$, in the sense of (\ref{dual}).

We consider $\lvec Y_t$ on an interval of $N$ sites with periodic boundary
conditions, started with an odd number of ones, and let $\lvec Y_\infty$
denote the process in equilibrium. For each $x\in\{0,1\}^N$, we define
\be\label{fxdef}
\rvec{f}_{x,N}(\al):=\frac{\P[|x\lvec Y_\infty|\mbox{ is odd}]}
{\P[\lvec Y_\infty(0)=1]}.
\ee
By our previous remarks, $\rvec{f}_{x,N}(\al)$, as a function of $x$ for fixed
$\al$, is a harmonic function for the one-sided rebellious voter model $\vec
X_t$, on an interval of $N$ sites with periodic boundary conditions. We will
be interested in the limit of $\rvec{f}_{x,N}(\al)$ as $N\to\infty$, when $x$
remains finite and fixed. In view of this, we only write down only the part of
$x$ that is nonzero. For example, we write $\rvec{f}_{1101,N}(\al)$ to denote
the function $\rvec{f}_{x,N}(\al)$ where $x$ is any element of $\{0,1\}^N$ such
that $(x(i),x(i+1),x(i+2),x(i+3))=(1,1,0,1)$ for some $i\in\{0,\ldots,N-1\}$
and $x(j)=0$ for all $j\neq i,i+1,i+2,i+3$.  Note that since the law of $\lvec
Y_\infty$ is translation invariant (modulo $N$), this definition does not
depend on the value of $i$.

Our simulations suggest that the functions $\rvec{f}_{x,N}(\al)$ converge as
$N\to\infty$ to a nontrivial limit function $\rvec{f}_x(\al)$, which, as a
function of $x$ for fixed $\al$, is a harmonic function for the one-sided
rebellious voter model $\rvec X_t$ on $\Z$, started with finitely many
ones. If $\al<\al_{\rm c}$, then on theoretical grounds one may expect that
\be\label{fro}
\rvec{f}_x(\al)=\frac{\rvec\rho_x(\al)}{\rvec\rho_1(\al)},
\quad\mbox{where}\quad
\rvec\rho_x(\al):=\P^x\big[\rvec X_t\neq\un 0\ \forall t\geq 0\big]
\ee
denotes the probability that the one-sided rebellious voter model $\rvec X_t$
started in the initial state $x$ survives. For $\al>\al_{\rm c}$, in the
regime where interface tightness holds, we can relate $\rvec{f}_x(\al)$ to the
invariant law of the process $\lvec Y_t$, started with a single one and viewed
from its left-most one (see formula (\ref{fromleft}) below).
Since we expect $\rvec{f}_x(\al)$ to be harmonic we expect that
$\rvec G\rvec{f}_x(\al)=0$, where $\rvec G$ is the generator of the one-sided
rebellious voter model, which is given by
\be
\rvec Gf_x=\sum_i\big(\al1_{\{x(i-1)\neq x(i)\}}
+(1-\al)1_{\{x(i-2)\neq x(i-1)\}}\big)\big(f_{x^{\{i\}}}-f_x\big).
\ee
In particular, the facts that $\rvec G\rvec{f}_1(\al)=0$ and $\rvec G\rvec
f_{11}(\al)=0$ lead to the relations
\bc
\dis\al\big(\rvec{f}_\emptyset(\al)-\rvec{f}_1(\al)\big)
+\big(\rvec{f}_{11}(\al)-\rvec{f}_1(\al)\big)
+(1-\al)\big(\rvec{f}_{101}(\al)-\rvec{f}_1(\al)\big)&=&0,\\[5pt]
\dis\al\big(\rvec{f}_{111}(\al)-\rvec{f}_{11}(\al)\big)
+\big(\rvec{f}_1(\al)-\rvec{f}_{11}(\al)\big)
+(1-\al)\big(\rvec{f}_{1101}(\al)-\rvec{f}_{11}(\al)\big)&=&0,
\ec
where $\rvec{f}_\emptyset(\al)=0$ and $\rvec{f}_1(\al)=1$. An inspection of our
numerical data shows that these equations are satisfied within the precision
of our simulations.

\begin{figure}
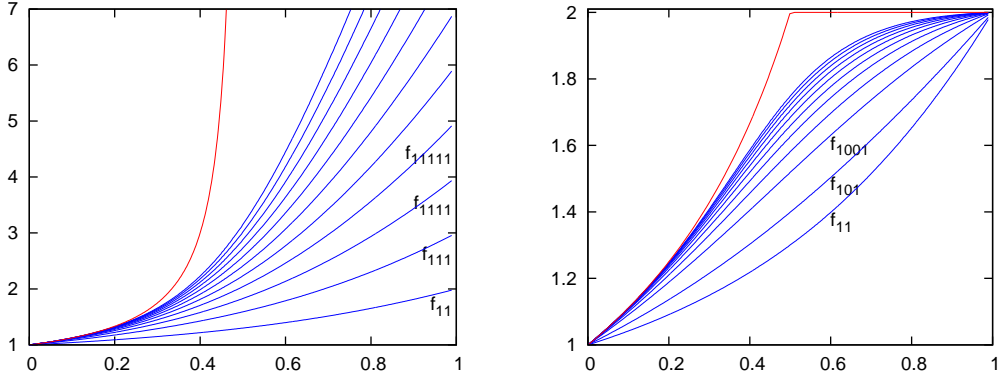

\begin{center}
\includegraphics[width=7cm]{\figfile{onemask2}}
\includegraphics[width=7cm]{\figfile{openmask}}
\caption
  [The functions $\vec{f}_{11},\vec{f}_{111},\vec{f}_{1111},\ldots$
  (left) and $\vec{f}_{11},\vec{f}_{101},\vec{f}_{1001},\ldots$ (right).]
  {The
  functions $\vec{f}_{11},\vec{f}_{111},\vec{f}_{1111},\ldots$ (left) and
  $\vec{f}_{11},\vec{f}_{101},\vec{f}_{1001},\ldots$ (right). The functions in
  the left figure increase to the limit $1/\rho(\al)$ while the functions in
  the right figure increase to $2-\rho(\al)$. Both limits (theoretical curves)
  are plotted together with the data. The functions
  $\vec{f}_{11},\vec{f}_{111},\vec{f}_{101}$
  are based on data obtained with the parameters $N=2^{11}$, $n=2^8$,
  $T=10^8$, $\al_{\rm b}=0$, $\al_{\rm e}=0.5$ (left half of the pictures) and
  $N=2^{15}$, $n=2^8$, $T=10^9$, $\al_{\rm b}=0.99$, $\al_{\rm e}=0.51$ (right
  half of the pictures). The other functions are a combination of simulations
  using the parameters $N=2^{11},n=2^8,T=10^8,\al_{\rm b}=0,\al_{\rm e}=0.5$
  and $N=2^{15},n=2^8,T=10^{10},\al_{\rm b}=0.99,\al_{\rm e}=0.51$.}
\label{fig:harmon}
\end{center}
\end{figure}

Our numerical data for the functions $\rvec{f}_x(\al)$, for some simple
choices of the pattern $x$, are shown in Figure~\ref{fig:harmon}. Note that
$\rvec{f}_1(x)$ (not shown in the figure) is by definition identically one,
because of the normalization chosen in (\ref{fxdef}). A surprising fact is
that the functions $\rvec{f}_x(\al)$ and their derivatives (shown in
Figure~\ref{fig:difharm}) appear to continue smoothly across the critical
point. Note that the normalization in (\ref{fxdef}) is crucial here, since
above $\al_{\rm c}$, both the numinator and the denuminator tend to zero as
$N\to\infty$. For $\al=1$ one has $\rvec{f}_x(\al)=|x|$, which is a well-known
harmonic function for the pure voter model. For $\al=0$, one has
$\rvec{f}_x(\al)=1$ for all nonzero $x$, reflecting the fact that the process
with $\al=0$ never dies out.

\begin{figure}
\begin{center}
\includegraphics[width=7cm]{\figfile{difmask4a}}
\includegraphics[width=7cm]{\figfile{difmask4b}}
\caption[Derivatives of the functions from
  Figure~\ref{fig:harmon}.]{Derivatives of the functions from
  Figure~\ref{fig:harmon}. On the left
  $\dif{\al}\rvec{f}_{11}(\al),\dif{\al}\rvec{f}_{111}(\al),\ldots$. On the
  right $\dif{\al}\rvec{f}_{11}(\al),\dif{\al}\rvec{f}_{101}(\al),
  \dif{\al}\rvec{f}_{1001}(\al),\ldots$. Figure based on the same data as
  Figure~\ref{fig:harmon}. The derivatives are estimated with a quadratic
  Savitzky-Golay filter using 11 data points.}
\label{fig:difharm}
\end{center}
\end{figure}

A few regularities may be observed from Figures~\ref{fig:harmon} and
\ref{fig:difharm}. Let 
\be
x_n:=\underbrace{111111111}_{n\mbox{ ones}}
\quad\mbox{and}\quad
x'_n:=1\hspace{-7pt}\underbrace{000000}_{n-2\mbox{ zeros}}\hspace{-7pt}1.
\ee
One expects that for $\al<\al_{\rm c}$,
\be
\lim_{n\to\infty}\rvec\rho_{x_n}(\al)=1
\quad\mbox{and}\quad
\lim_{n\to\infty}\rvec\rho_{x'_n}(\al)=1-(1-\rvec\rho_1(\al))^2.
\ee
Hence, in view of (\ref{fro}), one expects that
\be
\lim_{n\to\infty}\rvec{f}_{x_n}(\al)=\rvec\rho_1(\al)^{-1}
\quad\mbox{and}\quad
\lim_{n\to\infty}\rvec{f}_{x'_n}(\al)=2-\rvec\rho_1(\al),
\ee
which is indeed what we observe.

From Figure~\ref{fig:difharm} we moreover observe that
\be
\dif{\al}\rvec{f}_{x_n}(\al)\big|_{\al=1}=2(n-1)\qquad(n\geq 2).
\ee
This can be explained on theoretical grounds. Indeed, if we write
$\rvec{G}=\al\rvec{G}_{\rm vot}+(1-\al)\rvec{G}_{\rm reb}$, where
$\rvec{G}_{\rm vot}$ and $\rvec{G}_{\rm reb}$ are the (one-sided)
voter and rebellious parts of the generator, then a simple calculation shows
that
\be
\rvec{G}_{\rm vot}\dif{\al}\rvec{f}_x(\al)\big|_{\al=1}
=\,\rvec{G}_{\rm reb}\rvec{f}_x(1).
\ee
Setting $g_x:=\dif{\al}\rvec{f}_x(\al)\big|_{\al=1}$, this yields
\be\ba{l}\label{induc}
\dis g_{11}+g_\emptyset-2g_1=\,\rvec{G}_{\rm vot}g_1
=\,\rvec{G}_{\rm reb}\rvec{f}_1(1)=2,\\[5pt]
\dis g_{111}+g_1-2g_{11}=\,\rvec{G}_{\rm vot}g_{11}
=\,\rvec{G}_{\rm reb}\rvec{f}_{11}(1)=0,\\[5pt]
\dis g_{1111}+g_{11}-2g_{111}=\,\rvec{G}_{\rm vot}g_{111}
=\,\rvec{G}_{\rm reb}\rvec{f}_{111}(1)=0,\\[5pt]
\mbox{etcetera}
\ec
Since $\rvec{f}_\emptyset(\al)=0$ and $\rvec{f}_1(\al)=1$ we have
$g_\emptyset=g_1=0$ and we can solve $g_x$ inductively for
$x=11,111,1111,\ldots$ from (\ref{induc}).

\detail{
\be\ba{l}
\dis 0=\dif{\al}\big(\al\rvec{G}_{\rm vot}
+(1-\al)\rvec{G}_{\rm reb}\big)\rvec{f}_x(\al)\big|_{\al=1}\\[5pt]
\dis\quad=\big(\rvec{G}_{\rm vot}-\rvec{G}_{\rm reb}\big)
\rvec{f}_x(\al)\big|_{\al=1}
+\big(\al\rvec{G}_{\rm vot}+(1-\al)\rvec{G}_{\rm reb}\big)
\dif{\al}\rvec{f}_x(\al)\big|_{\al=1}\\[5pt]
\dis\quad=-\rvec{G}_{\rm reb}\rvec{f}_x(1)
+\rvec{G}_{\rm vot}\dif{\al}\rvec{f}_x(\al)\big|_{\al=1}.
\ec}

\begin{figure}
\begin{center}
\includegraphics[width=10cm]{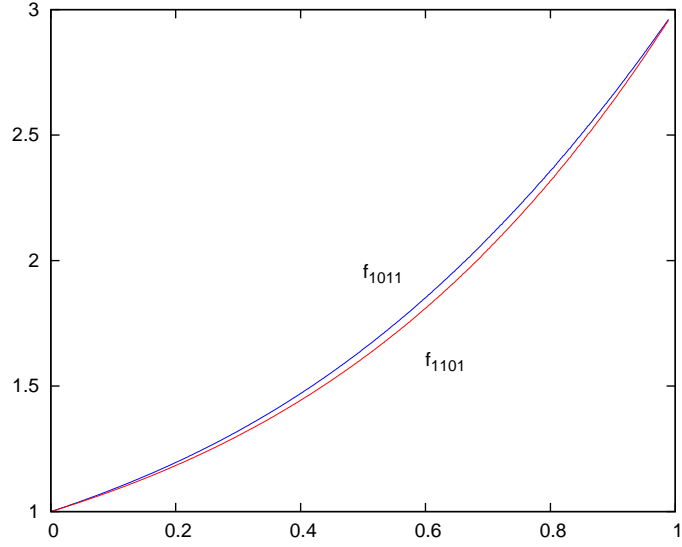}
\caption[The functions $\vec f_{1011}$ and $\vec f_{1101}$.]{The functions
  $\rvec{f}_{1011}$ and $\rvec{f}_{1101}$ for the one-sided model,
  demonstrating the asymmetry of $\rvec{Y}_\infty$. Combination of data
  obtained with the parameters $N=2^{11}$, $n=2^8$, $T=10^8$, $\al_{\rm b}=0$,
  $\al_{\rm e}=0.5$ and $N=2^{15}$, $n=2^8$, $T=10^9$, $\al_{\rm b}=0.99$,
  $\al_{\rm e}=0.51$.}
\label{fig:asym}
\end{center}
\end{figure}

The fact that the critical point and various other functions of the invariant
law are different for the one-sided and two-sided model implies that the law
of $\rvec Y_\infty$ cannot be mirror symmetric. Indeed, we observe that the
functions $\rvec{f}_{1011}(\al)$ and $\rvec{f}_{1101}(\al)$ are not identical, as show in
Figure~\ref{fig:asym}.

\subsection{Frequencies for three and more particles}\label{S:freq}

\begin{figure}
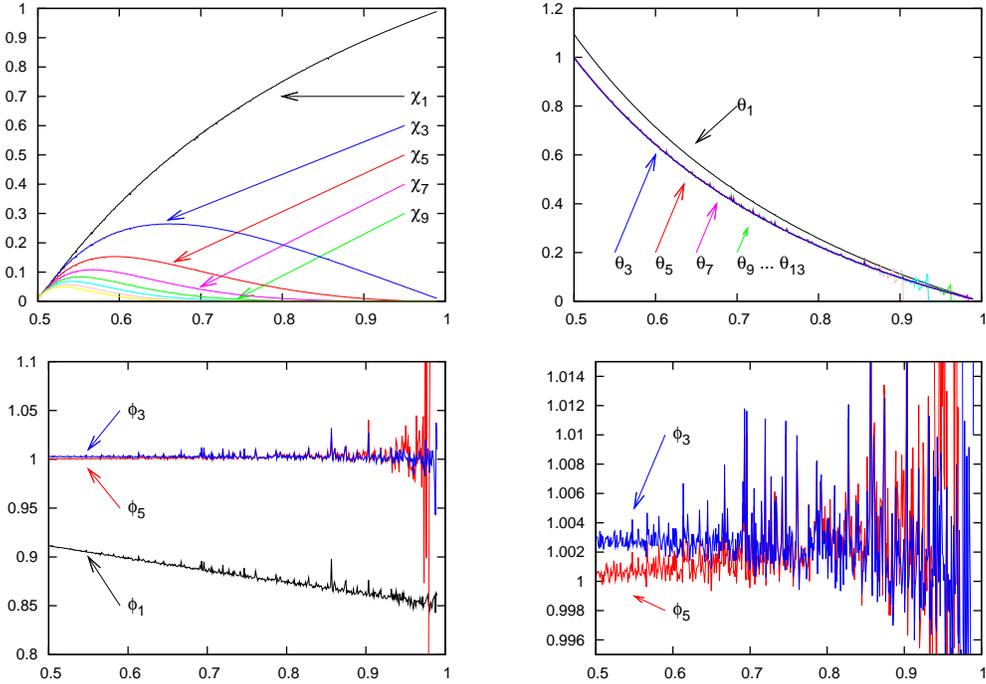

\begin{center}
\includegraphics[width=7cm]{\figfile{combi1}}
\includegraphics[width=7cm]{\figfile{combi2}}
\includegraphics[width=7cm]{\figfile{combi3}}
\includegraphics[width=7cm]{\figfile{combi4}}
\caption[Frequencies of one, three, five\ldots particles.]  {The frequency
  $\chi_k(\al)=\P[|Y_\infty|=k]$ of seeing $k$ particles for $k=1,3,5,\ldots$
  in the one-sided model. The second and third plot show the functions
  $\tet_k(\al):=\chi_{k+2}(\al)/\chi_k(\al)$ and
  $\phi_k(\al):=\tet_{k+2}(\al)/\tet_k(\al)$ for the first few values of
  $k=1,3,5,\ldots$. The fourth plot shows detail of the third plot. Data
  obtained with the parameters for $\chi$ listed in (\ref{sympar})}
\label{fig:threemore}
\end{center}
\end{figure}

As in the previous section, letting $\lvec Y_\infty$ denote the equilibrium
interface process for the one-sided rebellious voter model, we have plotted in
Figure~\ref{fig:threemore} (picture on top left) the functions
\be
\chi_k(\al)=\P\big[|\lvec Y_\infty|=k\big]\qquad(k=1,3,5,\ldots).
\ee
In particular, $\chi_1$ is the function $\chi$ from
(\ref{approx})~(ii). Except for the case $k=1$ (see (\ref{explicit})), we have
not found any simple explicit formulas that fit these curves. Nevertheless, some
regularities may be noted. In particular, if we define
\be\label{tetdef}
\tet_k(\al):=\frac{\chi_{k+2}(\al)}{\chi_k(\al)}
\quad\mbox{and}\quad
\phi_k(\al):=\frac{\tet_{k+2}(\al)}{\tet_k(\al)}\qquad(k=1,3,5,\ldots),
\ee
then it seems that the functions $\tet_k$ converge very fast to a limiting
function $\tet_\infty(\al):=\lim_{k\to\infty}\tet_k(\al)$. (See
Figure~\ref{fig:threemore} picture on top right; the functions
$\tet_3,\tet_5,\ldots$ all seem to fall on top of each other.) Likewise, the
functions $\phi_k$ converge very fast to the function that is identically
one. (See Figure~\ref{fig:threemore}, two bottom pictures.) It is hard to
obtain sufficiently precise data, but our best data (shown here) suggest that
the function $\phi_3$ is not identically one, although it is close. This means
that the functions $\tet_3$ and $\tet_5$ are probably not identical, even
though they almost fit each other in the picture on top right. Our numerical
data are well fitted by the linear/constant functions:
\bc\label{phiest}
\dis\phi_1(\al)&\approx&\dis c_0+c_1(\al-\ffrac{1}{2})\quad\mbox{with}\quad
c_0=0.9115\pm0.0015\quad\mbox{and}\quad c_1=-0.13\pm0.01,\\[5pt]
\dis\phi_3(\al)&\approx&\dis1.0027\pm0.0005.
\ec
and $\phi_k(\al)=1$ for $k\geq 5$. We note that if one assumes that the explicit
formula for $\chi_1(\al)$ in (\ref{explicit}) is correct, and one knows the
functions $\phi_k(\al)$ for each $k=1,3,\ldots$, then using the fact that
$\sum_{n=0}^\infty\chi_{2n+1}(\al)=1$ one has enough equations to solve the
functions $\chi_k(\al)$ for each $k=1,3,\ldots$.

The limiting function $\tet_\infty(\al):=\lim_{k\to\infty}\tet_k(\al)$ seems
to be be strictly decreasing and concave on $[\ffrac{1}{2},1]$ and satisfy
$\tet_\infty(\ffrac{1}{2})=1$, $\tet_\infty(1)=0$. This implies that the
distribution of the number of particles in equilibrium has an exponentially
decaying tail for each $\al>\frac{1}{2}$. In particular, the mean number of
particles
\be\label{mudef}
\mu(\al):=\E\big[|Y_\infty|\big]
\ee
is finite for all $\al>\frac{1}{2}$ and diverges as $\al\down\frac{1}{2}$.

The picture for the (two-sided) rebellious voter model is extremely similar,
except that the critical point is no longer $0.5$ and we have no
explicit formula for $\chi_1(\al)$.

\subsection{Edge speeds}\label{S:edge}

Let $Y_t$ be the interface model of the (two- or one-sided) rebellious voter
model started with a finite number of ones and let
\be
l_t:=\inf\{i\in\Z:Y_t(i)=1\}
\quad\mbox{and}\quad
r_t:=\sup\{i\in\Z:Y_t(i)=1\}
\ee
denote the position of the left-most and right-most ones, respectively.  Our
simulations suggest that there exist constants $v_-(\al)\leq v_+(\al)$, called
the {\em left} and {\em right edge speed}, respectively, such that
\be
\lim_{t\to\infty}\frac{l_t}{t}=v_-(\al)
\quad\mbox{and}\quad
\lim_{t\to\infty}\frac{r_t}{t}=v_+(\al)
\quad{\rm a.s.}
\ee
Set
\be\label{fromleft}
Z_t(i):=Y_t(l_t+i)\qquad(i\geq 0,\ t\geq 0).
\ee
Then $Z_t$ describes the process $Y_t$ as seen from the left-most one. It
seems reasonable to assume that the law of $Z_t$ converges as $t\to\infty$ to
some equilibrium law. Let $Z_\infty$ denote the process in equilibrium. The
left edge speed $v_-(\al)$ can be obtained as the equilibrium expectation
$\E[f(Z_\infty)]$ of a suitable function $f$ which sums the sizes of all
possible changes of $l_t$ multiplied with the rate at which they occur.

We have simulated the process $Z_t$ on a finite interval of $N$ sites by
neglecting all particles that fall off the right edge of the interval due to
the dynamics of $Y_t$ or the shifts in $l_t$. This process does not preserve
parity, but if $N$ is large enough the probability that all particles
annihilate each other is very small. For most of our data points, such events
never occurred, and when they occurred (especially near the critical point)
they were so rare that they likely had no big influence on our estimate of the
speeds. By slowly lowering or raising $\al$ in our usual fashion and keeping
track of all changes in $Z_t$ that correspond to a change in $l_t$ we obtained
numerical data for the left edge speed $v_-(\al)$ and, in a similar fashion,
also for the right edge speed $v_+(\al)$, both for the two-sided and one-sided
model.

\begin{figure}
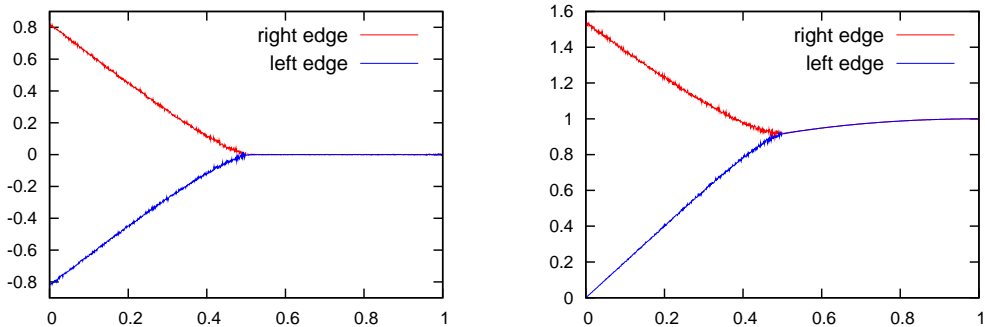

\begin{center}
\includegraphics[width=7cm]{\figfile{speed}}
\includegraphics[width=7cm]{\figfile{mspeed}}
\caption[Edge speeds for the rebellious voter model (left) and its one-sided
  counterpart (right).]{Edge speeds for the rebellious voter model (left) and
  its one-sided counterpart (right). The pictures show the speed of the right
  edge and left edge as a function of $\al$. Parameters used are
  $N=2^{11}$, $n=2^9$, $T=10^9$, $\al_{\rm b}=0$, $\al_{\rm e}=0.5$ (left half
  of each picture) and $N=2^{11}$, $n=2^9$, $T=10^{11}$, $\al_{\rm b}=0.999$,
  $\al_{\rm e}=0.5$ (right half of each picture).}
\label{fig:edge}
\end{center}
\end{figure}

\begin{figure}
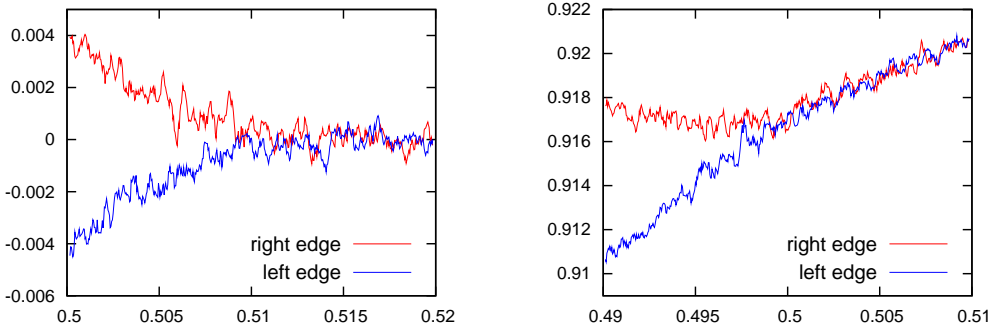

\begin{center}
\includegraphics[width=7cm]{\figfile{speedet_rev}}
\includegraphics[width=7cm]{\figfile{mspeedet_rev}}
\caption[Edge speeds; detail around critical point.]{Detail of
  Figure~\ref{fig:edge}. Parameters used are $N=2^{15}$, $n=2^9$, $T=10^{11}$,
  with $\al_{\rm b}=0.52$, $\al_{\rm e}=0.5$ for the two-sided model and
  $\al_{\rm b}=0.51$, $\al_{\rm e}=0.49$ for the one-sided model. The curves
  have been smoothed, corresponding to (effectively) $n=2^6$.}
\label{fig:edgedet}
\end{center}
\end{figure}

The results are shown in Figure~\ref{fig:edge}. One has $v_-(\al)<v_+(\al)$ if
and only if $\al<\al_{\rm c}$ and the functions $v_-(\al)$ and $v_+(\al)$ are
strictly increasing resp.\ decreasing on $[0,\al_{\rm c})$. For the two-sided
model, by symmetry, $v_-(\al)=-v_+(\al)$ and hence $v_-(\al)=0=v_+(\al)$ for
$\al\geq\al_{\rm c}$, while for the one-sided model one has $0\leq v_-(\al)$,
but otherwise the two pictures are remarkably similar. A detailed simulation
near the critical point (shown in Figure~\ref{fig:edgedet}) shows that the
estimates for the critical points which one obtains from these simulations are
consistent with our earlier estimates. Apart from the obvious values
$v_-(0)=0$ and $v_-(1)=v_+(1)=1$ (for the one-sided model) and $v_-=v_+=0$ on
$[\al_{\rm c},1]$ (for the two-sided model) we are not able to say anything
explicit about the curves $v_-(\al)$ and $v_+(\al)$.

\section{Conclusions}

We have simulated two one-dimensional models of parity preserving branching
and annihilating random walks, which arise as the interface models of the
(two-sided) rebellious voter model introduced in \cite{SS08hv} and the
one-sided rebellious voter model introduced in the present paper,
respectively. Both models appear to go through a phase transition between an
active and an inactive phase as the parameter $\al$ (which is one minus the
branching rate) is increased from zero to one. In most aspects, the processes
behave very similarly both qualitatively and quantitatively.

A peculiar property of the one-sided model, however, is that both the particle
density $\rho$ and the equilibrium probability of finding a single particle
$\chi$, defined in (\ref{rhodef}) and (\ref{chidef}), appear to be given by
explicit formulas, and the critical value appears to be exactly $\al_{\rm
  c}=\frac{1}{2}$. We have no idea why the one-sided model, which has less
symmetry than the two-sided model, should be more tractable than the
latter.\footnote{Perhaps the only aspect in which the one-sided model is
  potentially simpler than the two-sided model is that in the former,
  information is passed only from left to right. A third model, which we have
  only studied briefly, in which voter model updates look to the right and
  rebellious updates look to the left, seems to have a critical point
  $\al_{\rm c}=0.510\pm0.003$ and functions $\rho,\chi$ that are close, but
  not identical to those of the two-sided rebellious voter model.} In fact,
most hypothetical explanations that come to one's mind should apply to the
two-sided model as well, for which our numerical data convincingly show that
the formulas in (\ref{explicit}) do not fit. We have tried without success to
find explicit formulas for other functions of our processes.

An important question is whether our proposed formulas in (\ref{explicit}) for
the functions $\rho$ and $\chi$ of the one-sided model are exact or only good
approximations. In formulas (\ref{rhofrac}) and (\ref{phiest}), we have seen
examples of explicit formulas containing some `strange' constants that fit the
numerical data for certain functions of the one-sided and two-sided model,
respectively, within the precision of our simulations. In these examples, we
believe the given formulas are probably not exact but more or less
coincidentally close to the real functions. On the other hand, after trying
for a considerable while to find explicit formulas for other functions of the
processes, as we did, without success, one really starts to appreciate how
well the simple formulas in (\ref{explicit}) fit our data. The explicit
formula for $\rho$ has moreover been upheld in repeated more precise
simulations near the critical point. Thus, our present position is that these
formulas are probably exact, even though we have no theoretical explanation
for this.

An important consequence of our proposed formula for the particle density
$\rho$ of the one-sided rebellious voter model, if it is exact, is that it
implies that the order parameter critical exponent $\bet$ for this model is
one. By the principle of universlity, one then expects this to be true for the
whole PC universality class. It is presently generally believed that $\bet$ is
somewhat smaller than one. Recent estimates are $\bet=0.92(2)$ \cite{Hin00}
and $\bet=0.95(1)$ \cite{OM06}, although older numerical work in \cite{IT98}
yielded an estimate consistent with $\bet=1$. A thorough discussion of the
reliability of these estimates lies out of the scope of the present paper. We
stress that while our proposed formula for $\rho$ strongly suggests that
$\bet=1$, the most convincing data for the correctness of this formula are
obtained some distance away from the critical point. Our best data near the
critical point are not good enough to convincingly distinguish between
$\bet=0.92$ and $\bet=1$; see Figure~\ref{fig:logplot}.

It is currently believed that for the PC universality class (in dimension
one), $\beta/\nu_{\perp}=\frac{1}{2}$. If also $\bet=1$, this suggests that
the static critical exponents of this universality class are simple. (The
dynamical critical exponents, such as those related to the edge speeds, could
still be nontrivial.) In this context, it is interesting to note the effect,
observed in Section~\ref{S:harm}, that the harmonic function $\rvec{f}_x(\al)$
defined there appears to continue smoothly across the critical point, which
also suggests simple critical behavior.

We conclude this paper with some suggestions for further work. First of all,
the proposed formulas (\ref{explicit}) ask for a theoretical explanation. This
does not seem to be easy. The task would be easier if one could find explicit
formulas for more functions. In particular, if one could find explicit
formulas for the harmonic functions $\rvec{f}_x(\al)$ define in
Section~\ref{S:harm}, then it would be straightforward to prove that these are
indeed harmonic, which probably could be used to prove (\ref{explicit}) as
well. For those interested in such problems, we have made the data underlying
our figures available in the supplementary material. Even in the
absence of an explicit formula, these harmonic functions (and their observed
smooth behavior) seem a fruitful object for further theoretical considerations.

Our paper also calls for a critical evaluation of the evidence collected so
far about the order parameter critical exponent $\bet$ of the PC universality
class. Finally, it raises the question if there are more models in this
universality class that are (or appear to be) explicitly solvable. Our present
work suggests that such models might in particular be expected in the class of
one-sided models, where information is passed in one direction only.

\end{document}